\documentclass{article}
\usepackage{amsfonts,amssymb,amsmath}

\newtheorem{thm}{Theorem}
\newtheorem{dfn}{Definition}[subsection]
\newtheorem{prp}[dfn]{Proposition}
\newtheorem{lmm}[dfn]{Lemma}
\newtheorem{cor}[dfn]{Corollary}

\makeatletter

\@addtoreset{equation}{section}
\makeatother

\def \Z {{\mathbb Z}}

\newcommand{\cst}{{\rm cst\,}}
\newcommand{\one}{{\tt 1 \hspace{-0.8ex} \tt l}}
\newcommand{\qed}{\hspace*{2ex} \hfill $\Box$}

\newcommand{\Tau}{{\cal T}}

\title{An upper bound for front propagation velocities\\
inside moving populations}
\author{A.\ Gaudilli\`ere\footnote{
 Dipartimento di Matematica,
 Universit\`a di Roma Tre,
 Largo S.\ Leonardo Murialdo 1, 00146 Rome,
 Italy},\\
F.\ R.\ Nardi\footnote{
 EURANDOM, P.O.\ Box 513, 5600 MB Eindhoven,
 The Netherlands.
 }
 \footnote{
 Department of Mathematics and Computer Science,
 Eindhoven University of Technology,
 P.O.\ Box 513, 5600 MB Eindhoven,
 The Netherlands.
 }
}
\date{}

\begin{document}

\maketitle

\begin{abstract}

We consider a two type
(red and blue or $R$ and $B$)
particle population
that evolves on the $d$-dimensional
lattice according to some
reaction-diffusion process
$R+B\rightarrow 2R$
and starts with a single red particle
and a density $\rho$ of blue particles.
For two classes of models we give
an upper bound on the propagation velocity
of the red particles front
with explicit dependence on $\rho$.

In the first class of models
red and blue particles
respectively evolve
with a diffusion constant $D_R=1$
and a possibly time dependent jump rate
$D_B \geq 0$
-- more generally blue particles
follow some independent bistochastic process
and
this also includes long range random walks
with drift and various deterministic
processes.
We then get in all dimensions
an upper bound of order
$\max(\rho,\sqrt\rho)$
that depends only on $\rho$ and $d$
and not on the specific process
followed by blue particles,
in particular that does not
depend on $D_B$.
We argue that for $d \geq 2$
or $\rho \geq 1$
this bound can be optimal (in $\rho$),
while for the simplest case with $d=1$
and $\rho < 1$ known as the frog model,
we give a better
bound of order $\rho$.

In the second class of models
particles evolve with exclusion
and possibly attraction
inside a large two-dimensional box
with periodic boundary
conditions
according to Kawasaki dynamics
(that turns into simple exclusion 
when the attraction is set to zero.)
In a low density regime
we then get an upper bound
of order $\sqrt\rho$.
This proves a long-range decorrelation
of dynamical events
in this low density regime.

\end{abstract}

\section{Models and results}

\subsection{A diffusion-reaction model}

In \cite{KS} Kesten an Sidoravicius
considered the following Markov process.
A countable number of red and blue particles
perform independent continuous time
simple random walks
on the $d$-dimensional lattice $\Z^d$.
Red particles jump at rate $D_R$
and blue particles jump at rate $D_B$.
When a blue particle jumps on a site occupied
by a red particle, the blue particle
turns red. 
When a red particle 
jumps on a site occupied 
by blue particles these turn red.
Thinking respectively
at the red and blue particles
as individuals who have heard 
about a certain rumor and are ignorant
of it -- or as individuals who have
or have not a certain contagious disease --
this Markov process provides
a model of rumor propagation
-- or epidemic diffusion --
inside a moving population.
This is also a reaction-diffusion
dynamics of the kind $R+B\rightarrow 2R$
that can model a combustion process.

We define a each time $t\geq 0$ 
a {\em red zone} ${\cal R}(t)$
that is the set of sites $\Z^d$
that have been reached
by some red particle
at some time $s\in [0;t]$.
At any time $ t \geq 0$ 
all the red particles
stand in the red zone,
but some blue particles
can stand in the red zone
and the red zone can contain
empty sites.
The red zone is the set of sites
reached by the rumor
or the set of burnt sites according
to one or another interpretation
of the process.

Let us assume that the initial
configuration
was built in the following way.
We put independently in each site
$z\in\Z^d$
a random number of blue particles
according to Poisson variables
of mean $\rho>0$,
then at time $t=0$
we choose one particle 
according to some probabilistic
or deterministic rule,
we turn it red 
and we turn red the possible other
particles that stood in the same site.
Then, denoting by $B(z,r)$
the Euclidean ball of center $z$
and radius $r$ and making a change
of origin to have 
${\cal R}(0) = \{0\}$,
Kesten and Sidoravicius proved
\cite{KS}:

\medskip\par\noindent
{\bf Theorem [Kesten-Sidoravicius]:}
\em 
If $D_B=D_R>0$ there are two positive
and finite constants $C_1<C_2$
such that with probability 1
\begin{equation}
B(0,C_1 t) \subset {\cal R}(t)
\subset B(0,C_2 t)
\label{loup}
\end{equation}
will hold for all $t$
larger than some finite
random time $T_0$.

If $D_R>0$ there is a finite constant $C_2$
such that with probability 1
\begin{equation}
{\cal R}(t)
\subset B(0,C_2 t)
\label{up}
\end{equation}
will hold for all $t$
larger than some finite
random time $T_0$.
\rm

\medskip\par\noindent
{\bf Remarks:}
{\bf i)}
 Actually they did not introduced any change of origin.
The analogous result without change of origin
is an equivalent statement,
but our change of origin will serve us later.

\smallskip\par\noindent
{\bf ii)}
They proved the theorem in a slightly more general
situation: when the initial configuration
is obtained by adding any finite number of red particles
in a finite set set of sites
to a Poissonian distribution
of blue particles.
However it is easy to see
that the same result in this more general
case is equivalent to the previous theorem.
For the sake of simplicity
we will restrict ourselves
to discuss processes built like above,
starting with a single blue particle
that turns red.

\smallskip\par\noindent
{\bf iii)}
The inclusion (\ref{up})
gives a ``ballistic upper bound'' on ${\cal R}(t)$.
The ``ballistic lower bound'' expressed in (\ref{loup})
is much harder to prove and was obtain only in the special
cases $D_B=D_R>0$ \cite{KS} and $D_B=0$
(\cite{AMPR}, \cite{AMP}, \cite{RS}).
But it is believed that such a bound holds
in the general case $D_R>0$
(see \cite{Pa}).

\smallskip\par\noindent
{\bf iv)}
On the basis of (\ref{loup}), that is of a ballistic
upper {\em and} lower bound on ${\cal R}(t)$, Kesten
and Sidoravicius proved a ``shape theorem''
for the red zone: ${\cal R}(t)/t$
converges with probability 1
to a deterministic shape.
This proves the existence of a (maybe non isotropic)
propagation velocity of the rumor
or the combustion front.
In this context $C_1$ and $C_2$ are respectively
uniform lower and upper bounds
of this possibly non isotropic
front propagation velocity.

\smallskip\par\noindent
{\bf v)} 
It is believed that in the general case $D_R>0$
this propagation velocity does not depend on $D_B$
(see \cite{Pa} - note 38).

\medskip\par
In this paper we give an upper bound
on the propagation velocity,
i.e., a ballistic upper bound on 
${\cal R}(t)$ of the kind
(\ref{up}) with explicit
dependence of $C_2$ on the density
$\rho$ and no dependence on $D_B$.
This bound will be, in all dimensions,
of order $\max(\rho,\sqrt\rho)$.
We argue that for $d \geq 2$
or $\rho \geq 1$
this bound can be optimal (in $\rho$),
while for $d=1$
and $\rho < 1$,
we give in the simplest case
$D_B=0$ 
a better
bound of order $\rho$.
In addition we prove that
our upper bound in $\max(\rho,\sqrt\rho)$
holds for a larger class of models.
We prove it, on the one hand,
for those models in which red particles
perform independent random walks
while blue particles follow
any kind of independent bistochastic
process (see below).
On the other hand,
we give an analogous upper bound
for models in which the rumor diffuses
through a ``contact process''
inside an interacting particle system
with exclusion and possible attraction
(simple exclusion, Kawasaki dynamics)
when a low density limit
allows for a Quasi Random Walk approximation
as introduced in \cite{GdHNOS1}.

\subsection{One upper bound for many models}\label{thms}

We now define the first class of models
we will work with.
Like previously we start with
a density $\rho> 0$ of particles
putting independently in each site 
$z\in \Z^d$
a Poissonian number of particles
with mean $\rho$.
We then put labels 1, 2, 3, \dots\
on particles, we call $z_i$
the position of the particle $i$
and for all $t>0$
we will call $X_i(t)\in\Z^d$
and $Y_i(t)\in\{R;B\}$
its position and its color at time $t$.
With each $i$
we associate $Z_i^R$ and $Z_i^B$
two continuous time Markov processes
on $\Z^d$ in such a way that:
\begin{itemize}
\item all these processes start from 0
and are independent between them;
\item $Z_i^R$ is a simple random walk process
with diffusion constant
or jump rate~1;
\item $Z_i^B$ is a bistochastic process, i.e.,
satisfies
\begin{equation}
\forall z\in\Z^d, \forall t \geq 0,
\sum_{z_0\in\Z^d} P(z_0 + Z_i^B(t) = z) = 1
\label{bs}
\end{equation}
This includes simple random walks with constant 
or time dependent jumps rates, long range random walks
with drift, various deterministic processes,~\dots
\item
The $Z_i^B$'s (like the $Z_i^R$'s) have the same law.
\end{itemize}

At time $t=0$ we choose one particle
$i_0$ with some probabilistic or deterministic rule,
we change the origin to put it where $i_0$ stands,
we give the red color
to the particles in the new origin
and the blue color to the other particles
so that, for all $i$,
\begin{eqnarray}
X_i(0) &=& z_i - z_{i_0}\\
Y_i(0) &=& R\quad\mbox{if } X_i(0) = 0\\
Y_i(0) &=& B\quad\mbox{if } X_i(0) \neq 0
\end{eqnarray}

Then each particle $i$ follows the moves
of $Z_i^B$ while $Y_i=B$, turns red when it meets
a red particle and then follows the moves
of $Z_i^R$. More formally,
with for all $i$,
\begin{equation}
\tau_i := \left\{
  \begin{array}{l}
    0 \quad \mbox{if } X_i(0)=0\\
    \\
    \inf\{t\geq 0 :\: Y_i(t_-) = B,
    \exists j\neq i, Y_j(t_-) = R,  
    X_i(t) = X_j(t)\}\\
    \quad \mbox{if } X_i(0) \neq 0
  \end{array}
\right.
\end{equation}
with the usual convention $\inf\emptyset = +\infty$,
we have
\begin{eqnarray}
X_i(t) &=& \left\{
  \begin{array}{l}
    X_i(0) + Z_i^B(t) \mbox{ if } t\leq \tau_i\\
    X_i(0) + Z_i^B(\tau_i) + Z_i^R(t - \tau_i)
    \mbox{ if } t > \tau_i
  \end{array}
\right.\\
Y_i(t) &=& \left\{
  \begin{array}{l}
    B \mbox{ if } t < \tau_i\\
    R \mbox{ if } t \geq \tau_i
  \end{array}
\right.
\end{eqnarray}

We will call process of type RB any process
that can be built
in this way.
The Kesten and Sidoravicius reaction-diffusion model 
is a process of type RB when $D_R=1$.
We will call it KS process.
The general case $D_R>0$ can be mapped
on the KS process by a simple time rescaling.

Setting, like previously, for all $t\geq 0$,
\begin{equation}
{\cal R}(t) :=
\{z\in\Z^d :\: \exists i \geq 1, \exists s\in [0;t],
(X_i,Y_i)(s)=(z,R)\}
\end{equation}
we will prove

\begin{thm} \label{thm1}
There is a positive constant $\delta_d$ that depends
only on $d$ and such that,
for any RB process and for all $t\geq 0$
\begin{equation}
P\left(
  \exists z\in {\cal R}(t)\setminus B\left(
    0,\frac{\max(\rho,\sqrt\rho) t}{\delta_d}
  \right)
\right)
\leq
\frac{e^{-\delta_d\rho t}}{\delta_d\rho^4}
\end{equation}
\end{thm}

As a consequence, using the Borel-Cantelli lemma   
we get:

\begin{cor}
There is a positive constant $\delta_d$
that depends only on $d$
and such that for any RB process, with probability 1
\begin{equation}
{\cal R}(t)
\subset B\left(
  0,\frac{\max(\rho,\sqrt\rho) t}{\delta_d}
\right)
\end{equation}
will hold for all $t$ larger than some finite
random time $T_0$.
\end{cor}

We will give an analogous result
for a second class of models.
In dimension $d=2$ we consider
a low density lattice gas, with density
$\rho$, that evolves with exclusion and attraction
inside a large
finite box $\Lambda(\rho)$
with periodic boundary conditions
and according to the following Kawasaki dynamics
at inverse temperature
$\beta\geq 0$.
With
\begin{equation}
N:= \rho|\Lambda(\rho)|
\end{equation}
where $|\Lambda(\rho)|$
denotes the volume of $\Lambda(\rho)$,
we will write  $\hat\eta_i(t)\in\Lambda(\rho)$
for the position at time $t$ of the particle
$i$ in $\{1;\dots;N\}$
and $\eta_t\in\{0;1\}^{\Lambda(\rho)}$
for the configuration of the occupied sites
in $\Lambda(\rho)$,
in such a way that, for all $t\geq 0$,
\begin{equation}
\sum_{z\in\Lambda(\rho)}\eta_t(z) = N
\end{equation}
The energy of a configuration
$\eta\in\{0;1\}^{\Lambda(\rho)}$
is
\begin{equation}
H(\eta):=
\sum_{
  \stackrel{
    {\scriptstyle \{x;y\}\in\Lambda(\rho)}
  }{
    {\scriptstyle |x-y|=1}
  }
}
-U\eta(x)\eta(y)
\end{equation}
where $|\cdot|$ stands now for the Euclidean norm
and $-U\leq 0$ is the binding energy.
With each particle we associate a Poissonian clock
of intensity 1. At each time $t$ when a particle's
clock rings we extract with uniform
probability a nearest neighbor site of the particle, 
say $i$.
If this site is occupied by another particle
then $i$ does not move.
If not, we consider the configuration
$\eta'$
obtained by moving $i$ to the vacant site,
then with probability 
\begin{equation}
p=e^{-\beta[H(\eta')-H(\eta)]_+}
\end{equation}
$i$ moves to the vacant site
and, with probability $1-p$,
$i$ remains where it was at time $t_-$.
Observe that the case $U=0$
corresponds to the simple exclusion process.

In addition we choose at time $t=0$
some particle $i_0$
according to some probabilistic
or deterministic rule and give
to $i_0$, as well as to the particles
that share with $i_0$ the same cluster
at time $t=0$,
the red color,
while all the other particles
receive
the blue color.
Like previously a red particle
will definitively remain red
and a blue particle turns red
as soon as it shares some cluster
with some red particle.
We call RBK process
this dynamics
and, for all $t\geq 0$,
the red zone ${\cal R}(t)$
is defined like above.

To control the propagation of the red particles
in the regime $\rho\rightarrow 0$
we will use the low density to reduce
the problem to simple random walks estimates.
This is more challenging
when $\rho$ and $\beta$ go jointly
to 0 and $+\infty$: in this case we have not
only a low density regime but also
a strong interaction regime.
We will then deal with this more
challenging regime only,
setting $\rho=e^{-\Delta\beta}$
for $\Delta$ a positive parameter
and sending $\beta$ to infinity.
We will write $\Lambda_\beta$
for $\Lambda(\rho)$ and we will choose
$|\Lambda_\beta|=e^{\Theta\beta}$
for some real parameter 
$\Theta>\Delta$.
This regime was studied in \cite{GdHNOS1}
where a ``Quasi Random Walk (QRW) property''
was proved up to the first time
of ``anomalous concentration'' 
$\Tau_{\alpha,\lambda}$.
For $\alpha$ a positive parameter
that can be chosen as close as 0
as we want and $\lambda$ a slowly
increasing and unbounded function
such that 
\begin{equation}
\lambda(\beta)\ln\lambda(\beta) = o(\ln\beta)
\label{clambda}
\end{equation}
(for example $\lambda(\beta)=\sqrt{\ln\beta}$),
$\Tau_{\alpha,\lambda}$
is defined as the first time
when there is a square box
$\Lambda\subset\Lambda_\beta$
with volume less than $e^{\beta(\Delta-\alpha/4)}$
that contains more than $\lambda/4$ particles.
We will recall and use this QRW property
to prove

\begin{thm} \label{thm2}
For the RBK process,
for all $\delta>0$
and all $C>0$,
uniformly in the starting configuration,
and uniformly in $T=T(\beta)\leq e^{C\beta}$
\begin{equation}
P\left(
  \exists z\in {\cal R}(T)\setminus B\left(
    0,e^{\delta\beta}\sqrt{\rho}T
  \right)
  \mbox{ and }
  \Tau_{\alpha,\lambda}>T
\right)
\leq
\rho^{-3}e^{\delta\beta}\exp\left\{
  -e^{-\delta\beta}\rho t
\right\}
+SES
\end{equation}
where $SES$ stands for ``super exponentially small'',
i.e., for
a positive function $f$ that does not
depend on $T$ and the starting configuration
and such that
\begin{equation}
\lim_{\beta\rightarrow +\infty}
\frac{1}{\beta}
\ln f(\beta) = - \infty
\end{equation}
\end{thm}

As a straightforward consequence:

\begin{cor}
\label{kcor}
For the RBK process,
for all $\delta>0$
and all $C_2>C_1>\Delta$,
uniformly in the starting configuration,
and uniformly in $T=T(\beta)$
such that $T\geq e^{C_1\beta}$
and $T\leq e^{C_2\beta}$,
\begin{equation}
P\left(
  \exists z\in {\cal R}(T)\setminus B\left(
    0,e^{\delta\beta}\sqrt{\rho}T
  \right)
\mbox{ and }
\Tau_{\alpha,\lambda}>T\right)
\leq
SES
\end{equation}
\end{cor}

Of course these results would
be of no use if we were not able
to have some control on $\Tau_{\alpha,\lambda}$.
But in \cite{GdHNOS1} we discussed
the fact starting from a ``good configuration''
$\Tau_{\alpha,\lambda}$ is ``very long''.
For example we proved that in the case $\Delta>2U$,
starting
from the canonical Gibbs measure
associated with $H$,
for all $C>0$,
\begin{equation}
P\left(
  \Tau_{\alpha,\lambda} < e^{C\beta}
\right)
=SES
\end{equation}

As a consequence of these results
we will prove a long range decorrelation
of dynamical events in this low density
regime.
Given $\Lambda^{(1)}$ and $\Lambda^{(2)}$
two square boxes contained in $\Lambda_\beta$
we will denote by
$d(\Lambda^{(1)},\Lambda^{(2)})$
their Euclidean distance
and by $({\cal F}_t^{(1)})_{t\geq 0}$
and $({\cal F}_t^{(2)})_{t\geq 0}$
the filtrations generated by 
$(\eta_{t\wedge \Tau_{\alpha,\lambda}}|_{
\Lambda^{(1)}})_{t\geq 0}$
and
$(\eta_{t\wedge \Tau_{\alpha,\lambda}}|_{
\Lambda^{(2)}})_{t\geq 0}$.
With these notations:

\begin{thm}\label{thm3}
For the Kawasaki dynamics,
for all $\delta>0$ and all $C>0$,
uniformly in the starting configuration,
uniformly in $T=T(\beta)\leq e^{C\beta}$,
uniformly in $\Lambda^{(1)}$ and $\Lambda^{(2)}$
such that
\begin{equation}
d(\Lambda^{(1)},\Lambda^{(2)})
\geq e^{\delta\beta}\max(\sqrt T, \sqrt\rho T)
\end{equation}
and uniformly
in $(A^{(1)},A^{(2)})\in
{\cal F}_T^{(1)}\times{\cal F}_T^{(2)}$
\begin{equation}
\left|
  P\left(A^{(1)}\cap A^{(2)}\right)
  -
  P\left(A^{(1)}\right)P\left(A^{(2)}\right)
\right|
\leq
SES
\end{equation}
\end{thm}

In the study of the low temperature
metastable Kawasaki dynamics
(the case $U<\Delta<2U$,
see \cite{dHOS})
we will need such a long range decorrelation
property (see \cite{GdHNOS1}).
This constituted the initial motivation
of this paper.

\subsection{How good are our bounds?}

In this paper we will not give
any lower bound on the propagation velocity.
But we give here some heuristic that indicates
that $\max(\rho,\sqrt\rho)$
should be the right order of the velocity
propagation in different situations.

Consider for now the KS process in dimension
$d=2$ with $\rho< 1$ and in the special case
$D_B=D_R=1$.
${\cal R}(t)$ should then look like a kind
of ball that contains all the red particles
and very few blue particles.
In addition $D_B=D_R$ implies that,
except for the color propagation,
the particle system starts
and remains at equilibrium.
Let us call $n(t)$ the number of red particles
at time $t$.
Since only the particles at the border
of ${\cal R}(t)$ should contribute to the propagation
of the rumor and since a particle typically
waits for a time $1/\rho$ before meeting
another particle, we should have
\begin{equation}
dn \simeq \cst\sqrt n \rho dt 
\label{heur}
\end{equation}
where `$\cst$' stands for a positive constant
the value of which can change from line to line.
As a consequence 
\begin{equation}
\sqrt n \simeq \cst\rho t 
\end{equation}
If $r(t)$ stand for the radius of the smallest
Euclidean ball that contains ${\cal R}(t)$
we should have
\begin{equation}
n \simeq \cst r^2 \rho  
\label{heurr}
\end{equation}
so that
\begin{equation}
r\simeq \cst\frac{\sqrt n}{\sqrt\rho}
\simeq \cst\sqrt\rho t 
\end{equation}

If $\rho\geq 1$ we will typically have $\rho$
particles per site and (\ref{heur}) turns into
\begin{equation}
dn \simeq \cst\rho\sqrt{\frac{n}{\rho}}\rho dt 
\end{equation}
so that
\begin{equation}
r\simeq \cst\sqrt{\frac{n}{\rho}}\simeq \cst\rho t 
\end{equation}

If $d\geq 3$ or $D_R\neq D_B$ we do not have
such kind of heuristic. 
In the former case indeed ${\cal R}(t)$ 
should be a more complex fractal object,
in the latter case the system does not stay
at equilibrium.
However Theorem \ref{thm1}
says that an upper bound of order 
$\max(\rho,\sqrt\rho)$ holds independently of $D_B$
and the dimension.

For $d=1$, $D_R=D_B$
and $\rho<1$ the previous heuristic
has to be modified.
In this case the typical inter-particle distance
is $1/\rho$ and a particle typically waits for a time
$1/\rho^2$ before meeting another particle.
Then (\ref{heur}) and (\ref{heurr}) turn into
\begin{eqnarray}
dn &\simeq & \cst \rho^2 dt \\
n &\simeq & \cst r \rho  
\end{eqnarray}
and we get
\begin{equation}
r\simeq \cst\rho t 
\end{equation}
while Theorem \ref{thm1}
gives only an upper bound
on the velocity of order
$\sqrt\rho > \rho$.
We will prove an upper bound
of order $\rho$
for the simplest case of the KS process,
that is $D_B=0$,
also known as {\em frog model}.

\begin{prp}
For the KS process in dimension 1 and with $D_B=0$,
there is a positive constant $\delta$
such that for all $t\geq 0$
\begin{equation}
P\left(
  \exists z\in {\cal R}(t)\setminus B\left(
    0,\frac{\rho t}{\delta}
  \right)
\right)
\leq
\frac{e^{-\delta\rho^2 t}}{\delta\rho^2}
\end{equation}
\label{pf1}
\end{prp}

As previously we then get
with the Borel-Cantelli lemma:

\begin{cor}
For the KS process in dimension 1
and with $D_B=0$
there is a positive constant $\delta$
such that, with probability 1
\begin{equation}
{\cal R}(t)
\subset B\left(
  0,\frac{\rho t}{\delta}
\right)
\end{equation}
will hold for all $t$ larger than some finite
random time $T_0$.
\end{cor}

We will give
in section \ref{gc}
some indications
on how one can extend the simple
proof of Proposition \ref{pf1}
to the general case of the KS processes.
This is rather technical
and we will not go beyond these indications.

\subsection{Notation and outline of the paper}

We will write $\cst$
for a finite and positive constant that 
depends only on the dimension $d$
and the value of which
can change from line to line.
Given $d\geq 1$ we will write $|\cdot |$
for the $d$-dimensional Euclidean norm.
Given a Markov process $X$
and $x$ in its state space,
we will write $P_x$
for the law of the process
that starts from $x$.

In section \ref{prlm} we prove simple random walk
and large deviations estimates and we recall
some definitions and properties regarding the QRW
approximation for the Kawasaki dynamics.
In section \ref{frog}
we prove Theorem \ref{thm1}
for the frog model
as well as Proposition \ref{pf1}.
In section \ref{gc} we prove Theorem \ref{thm1}
in the general case
as well as Theorems \ref{thm2} and \ref{thm3}.

\section{Preliminaries}
\label{prlm}

\subsection{Random walk and large deviation estimates}

\begin{lmm}
Let $N$ and $N'$ be two independent Poisson variables
and $\gamma >1$ such that $E[N']\geq \gamma E[N]$.
Then
\begin{eqnarray}
i) && P(N\geq \gamma E[N])
\leq \exp\{-E[N](\gamma\ln\gamma - (\gamma -1))\}
\label{ldi}\\
ii) && P\left(
  N \leq \frac{E[N]}{\gamma}
\right)
\leq \exp\left\{
  -E[N]\left(
    \left(
        1-\frac{1}{\gamma}
    \right) 
    -\frac{\ln\gamma}{\gamma}
  \right)
\right\}
\label{ldii}\\
iii) && P\left(
  \frac{N}{E[N]} \geq \gamma\frac{N'}{E[N']}
\right)
\leq 2\exp\left\{
  -E[N](t\ln t - (t-1))
\right\}
\label{ldiii}\\
&& \mbox{with } t:=\frac{\gamma-1}{\ln\gamma}
\in ]1;\gamma[
\nonumber
\end{eqnarray}
\end{lmm}

\noindent
{\bf Proof:} We just use
the Chebyshev exponential inequality.
With $\lambda=E[N]$ we have, for any $t\geq 0$,
\begin{equation}
P(N\geq \gamma\lambda)
\leq
e^{-t\gamma\lambda}E[e^{tN}]
=
\exp\{-\lambda(t\gamma-(e^t-1))\}
\end{equation}
Optimizing in $t$ we find (\ref{ldi})
with $t=\ln\gamma$.
Similarly, for any $t\geq 0$,
\begin{equation}
P(N\leq \lambda/\gamma)
\leq
e^{t\lambda/\gamma}E[e^{-tN}]
=
\exp\{-\lambda((1-e^{-t})-t/\gamma)\}
\end{equation}
Optimizing in $t$ we find (\ref{ldii})
with $t=\ln\gamma$.
Finally we have, for any $t\geq 0$,
\begin{equation}
P\left(
  \frac{N}{E[N]} \geq \gamma\frac{N'}{E[N']}
\right)
\leq 
P\left(
  N \geq tE[N]
\right)
+
P\left(
  N' \leq \frac{t}{\gamma}E[N']
\right)
\end{equation}
By (\ref{ldi}) and (\ref{ldii})
this gives, if $t>1$ and $t<\gamma$,
\begin{eqnarray}
&&P\left(
  \frac{N}{E[N]} \geq \gamma\frac{N'}{E[N']}
\right)
\nonumber\\
&&\quad\leq\quad 
\exp\{-\lambda(t\ln t - (t-1))\}
\nonumber\\
&&\quad\quad\quad
+\exp\left\{
  -\lambda\gamma\left(
    \left(
        1-\frac{t}{\gamma}
    \right) 
    +\frac{t}{\gamma}\ln\frac{t}{\gamma}
  \right)
\right\}
\end{eqnarray}
The two terms of this sum are equal
when 
\begin{equation}
t=\frac{\gamma-1}{\ln\gamma}
\end{equation}
The concavity of the logarithm ensures
\begin{equation}
1-\frac{1}{\gamma}
\leq -\ln\frac{1}{\gamma}
= \ln\gamma
\leq \gamma-1
\end{equation}
so that $1<t<\gamma$
and this gives (\ref{ldiii}).
\qed

\begin{lmm}\label{dbb}
Let $\zeta$ be a $d$-dimensional continuous
time simple random walk
with rate jump 1.
For all $t\geq 0$
and $z\in\Z^d$
\begin{itemize}
\item
$\mbox{if $|z| \leq t$ then}$
\begin{equation}
P_0(\zeta(t) = z)
\leq 
\frac{\cst}{t^{d/2}}\exp\left\{
  -\frac{\cst|z|^2}{t}
\right\}
\label{db}
\end{equation}
\item
$\mbox{if $|z| \geq t$ then}$
\begin{equation} 
P_0(\zeta(t) = z)
\leq 
{\cst}\exp\left\{
  -\cst|z|
\right\}
\label{bb}
\end{equation}
\end{itemize}
\end{lmm}

\noindent
{\bf Remark:}
Since we just need an upper bound
on these probabilities
we do not need the usual condition
$|z|=o(t^{2/3})$
of the local central limit theorem.
However, working with continuous time random walks,
we have to treat separately
the case $|z|>t$.

\medskip\par\noindent
{\bf Proof of the lemma:}
We will prove slightly different
but equivalent estimates:
(\ref{db}) when $|z| \leq 2t$ 
(\ref{bb}) when $|z| \geq 2t$.

For the case $|z| \geq 2t$ we apply
the previous lemma. If $\zeta$
reaches $z$ in time $t$
then the number of its clock rings
up to time $t$ is larger than
or equal to $|z|$.
Since this number has a Poissonian
distribution of mean $t$,
this occurs, by (\ref{ldi}),
with a probability smaller than
\begin{equation}
\exp\left\{
  -t\left(
    \frac{|z|}{t}\ln\frac{|z|}{t} - \left(
      \frac{|z|}{t}-1
    \right)
  \right)
\right\}
\leq 
\exp\left\{
  -t\frac{\cst|z|}{t}
\right\}
=
e^{-\cst|z|}
\end{equation}
(for the last inequality 
we used that $|z|/t$
was bounded away from 1.)

For the case $|z| \leq 2t$
we first observe that,
working with a continuous time
process with independent coordinates,
it is enough to prove the result
for $d=1$.
Then we prove the estimate for
$\tilde\zeta$
the discrete time version
of such a one dimensional process.
Without loss of generality
we can assume that $z\in\Z$
is non negative.
If $z\leq {n}/{2}$,
then, by the Stirling formula,
\begin{eqnarray}
&&P_0\left(
  \tilde\zeta(n) = z
\right)
\nonumber\\
&&\quad\leq\quad \frac{\cst}{2^n}
\frac{n^{n}e^{-n}\sqrt n}
{\left(\frac{n+z}{2}\right)^{\frac{n+z}{2}}
e^{-\frac{n+z}{2}}\sqrt{\frac{n+z}{2}}
\left(\frac{n-z}{2}\right)^{\frac{n-z}{2}}
e^{-\frac{n-z}{2}}\sqrt{\frac{n-z}{2}}}\\
&&\quad\leq\quad
\frac{\cst}{\sqrt n}
\frac{2}
{\sqrt{1+\frac{z}{n}}\sqrt{1-\frac{z}{n}}}
\left[
  \left(1+\frac{z}{n}\right)^{\frac{1+z/n}{2}}
  \left(1-\frac{z}{n}\right)^{\frac{1-z/n}{2}}
\right]^{-n}\\
&&\quad\leq\quad
\frac{\cst}{\sqrt n}
\exp\left\{-nI(z/n)\right\}
\end{eqnarray}
with
\begin{equation}
I(x) :=
\frac{1+x}{2}\ln(1+x)
+
\frac{1-x}{2}\ln(1-x),
\quad x\in[-1;1]
\end{equation}
It is immediate to check
that 
\begin{equation}
\left\{
  \begin{array}{l}
    I(0) = I'(0) = 0\\
    \forall x \in ]-1;1[, I''(x)=\frac{1}{1-x^2}\geq 1
  \end{array}
\right.
\end{equation}
As a consequence, for all $x\in[-1;1]$,
\begin{equation}
I(x)\geq \frac{x^2}{2}
\end{equation}
and this gives, for $z\leq {n}/{2}$,
\begin{equation}
P_0(\tilde\zeta(n) = z)
\leq 
\frac{\cst}{\sqrt n}\exp\left\{
  -\frac{z^2}{2n}
\right\}
\end{equation}
This is easily extended
to the case $z\geq n/2$, i.e.,
$z/n\geq 1/2$:
\begin{eqnarray}
P_0\left(
  \tilde\zeta(n) = z
\right)
&\leq & 
\cst\exp\left\{-nI(z/n)\right\}\\
&\leq & 
\cst\exp\left\{-n\cdot 8I(1/2)\frac{z^2}{2n^2}\right\}\\
&\leq & 
\cst\sqrt{\frac{n}{z^2}}
\exp\left\{-\frac{z^2}{2n}\right\}\\
&\leq & 
\frac{\cst}{\sqrt n}
\exp\left\{-\frac{z^2}{2n}\right\}
\end{eqnarray}
Finally we use the previous lemma
to prove (\ref{db}).
We have
\begin{eqnarray}
P_0\left(
  \zeta(n) = z
\right)
&\leq & 
E\left[
  \frac{\cst}{\sqrt N}
  \exp\left\{-\frac{z^2}{2N}\right\}
\right]
\end{eqnarray}
where $N$ is a Poissonian variable of mean $t$.
By (\ref{ldi}), (\ref{ldii})
applied with a large enough $\gamma$
we can find two positive constants $c_1$,
$c_2$ with $4c_1<c_2$
such that
\begin{eqnarray}
P_0\left(
  \zeta(n) = z
\right)
&\leq & 
  \frac{\cst}{\sqrt t}
  \exp\left\{-c_1\frac{z^2}{t}\right\}
  + \exp\left\{-2c_2t\right\}\\
&\leq & 
\frac{\cst}{\sqrt t}\left(
  \exp\left\{-c_1\frac{z^2}{t}\right\}
  + \exp\left\{-c_2t\right\}
\right)
\end{eqnarray}
and we get (\ref{db}) using 
$z\leq 2t$, i.e., $4t\geq z^2/t$.
\qed

\subsection{Quasi Random Walks} \label{qrw}

With the notation we introduced in section
\ref{thms} for the Kawasaki dynamics
and given an arbitrarily small parameter $\alpha>0$
as well as an unbounded slowly increasing 
function $\lambda$ satisfying (\ref{clambda}),
we recall in this section a few definitions
and results from
\cite{GdHNOS1}.

\begin{dfn}
A process $Z=(Z_1;\dots;Z_N)$ on $\Lambda_\beta^N$
is called a random walk with pauses (RWP) associated
with the stopping times
\begin{equation}
0=\sigma_{i,0} = \tau_{i,0}
\leq \sigma_{i,1} \leq \tau_{i,1}
\leq \sigma_{i,2} \leq \tau_{i,2}
\leq \dots
\quad
i\in\{1,\dots,N\}
\end{equation}
if for any $i$ in $\{1;\dots;N\}$,
$Z_i$ is constant on all time intervals
$[\sigma_{i,k},\tau_{i,k}]$, $k\geq 0$,
and if the process 
$\tilde Z=(\tilde Z_1,\dots,\tilde Z_N)$
obtained from $Z$;
by cutting off these pauses intervals,
i.e., with 
\begin{equation}
\tilde Z_i(s) := Z_i\left(
  s+\sum_{k<j_i(s)} \tau_{i,k}-\sigma_{i,k}
\right),
\quad s\geq 0
\end{equation}
where
\begin{equation}
j_i(s):=
\inf\left\{
  j \geq 0 :\:
  s+\sum_{k<j} \tau_{i,k}-\sigma_{i,k}
  \leq \sigma_{i,j}
\right\}
\end{equation}
is an independent random walk process in law.
\end{dfn}

Now with
\begin{equation}
T_\alpha:= e^{(\Delta-\alpha)\beta}
\end{equation}
Quasi Random Walk processes are defined
as follows.

\begin{dfn}
We say that a process $\xi=(\xi_1,\dots,\xi_N)$
on $\Lambda_\beta^N$
is a Quasi Random Walk process
with parameter $\alpha>0$ up to
a stopping time $\Tau$,
written QRW($\alpha$,$\Tau$),
if there exists a coupling between $\xi$
and a RWP process $Z$
associated with stopping times
\begin{equation}
0=\sigma_{i,0} = \tau_{i,0}
\leq \sigma_{i,1} \leq \tau_{i,1}
\leq \sigma_{i,2} \leq \tau_{i,2}
\leq \dots
\quad
i\in\{1;\dots;N\}
\end{equation}
such that $\xi(0)=Z(0)$,
for any $i$ in $\{1,\dots,N\}$
$\xi_i$ and $Z_i$ evolves jointly
($\xi_i-Z_i$ is constant)
outside the pause intervals
$[\sigma_{i,k},\tau_{i,k}]$,
$k\geq 0$, and for any $t_0\geq 0$
the following events occur with
probability $1-SES$ uniformly
in $i$ and $t_0$:
\begin{eqnarray}
F_i(t_0) &:=&
\Big\{
  \sharp\left\{
    k\geq 0 :\:
    \sigma_{i,k}\in\left[
      t_0\wedge\Tau, (t_0+T_\alpha)\wedge\Tau
    \right]
  \right\}
  \leq l(\beta)
\Big\}\\
G_i(t_0) &:=&
\Big\{
  \forall k\geq 0, \forall t\geq t_0,
  \sigma_{i,k}\in\left[
      t_0\wedge\Tau, (t_0+T_\alpha)\wedge\Tau
    \right]
\nonumber\\
&&\qquad  
  \Rightarrow
  \left|
    \xi(t\wedge\tau_{i,k}\wedge\tau)
    _ \xi(t\wedge\sigma_{i,k}\wedge\tau)
  \right|
  \leq l(\beta)
\Big\}
\end{eqnarray}
for some $\beta\mapsto l(\beta)$ that satisfies
\begin{equation}
\lim_{\beta\rightarrow +\infty}
\frac{1}{\beta}\ln l(\beta)=0
\end{equation}
\end{dfn}

In words, the fact that 
for each $i$ the events
$F_i(t_0)$ and $G_i(t_0)$
occur for all $t_0\geq 0$
means, on the one hand, that in each time interval
before time $\Tau$
and of length $1/\rho$ almost,
there are few pauses for the associated RWP $Z_i$
(a non exponentially large number)
and, on the other hand,
that $\xi_i$ stays close to $Z_i$
in the sense that during each
of these few pause intervals
the distance between the two processes
cannot increase of more than the same
non exponentially large quantity $l$.

\begin{prp}
For any unbounded and slowly increasing
function $\lambda$ that satisfies
(\ref{clambda}) and any positive
$\alpha<\Delta$,
$\hat\eta$ is a QRW($\alpha$,$\Tau_{\alpha,\lambda}$)
process.
\end{prp}

We refer to \cite{GdHNOS1} for the proof.
As a consequence of this QRW property
we have for all $\delta>0$, uniformly
in the initial configuration
and uniformly in $T=T(\beta)\in [2,T_\alpha^2]$
\begin{equation}
P\left(
  \Tau_{\alpha,\lambda}>T,
  \exists t\in [0,T], \exists i \in \{1,\dots,N\},
  |\hat\eta_i(t)-\hat\eta_i(0)|>e^{\delta\beta}T
\right)
\leq SES
\label{nnsupdiff}
\end{equation}
In \cite{GdHNOS1} we also introduced at any time
$t_0\geq 0$
a partition of $\{1;\dots;N\}$
in clouds of potentially interacting particles
on time scale $T_\alpha$:
we associate with each particle $i$
a ball centered at its position at time $t_0$
with radius 
\begin{equation}
r:=e^{\frac{\alpha}{4}\beta}\sqrt{T_\alpha}
\end{equation}
we call $B_0$
their union
\begin{equation}
B_0 := \cup_i B(\hat\eta_i(t_0),r)
\end{equation}
and we say that two particles are in the same cloud
if there are, at time $t_0$,
in the same connected component
of $B_0$.
It is easy to check that if
$t_0<\Tau_{\alpha,\lambda}$
then no cloud contains more than $\lambda$
particles.
And, as a consequence of (\ref{nnsupdiff}),
with probability $1-SES$
interactions between particles
during the time interval
$[t_0,(t_0+T_\alpha)\wedge\Tau_{\alpha,\lambda}[$
will only take place {\em inside} the different
clouds (and not between particles of different clouds.)

\section{The frog model}\label{frog}

\subsection{Proof of Theorem \ref{thm1}
for the KS process with $D_B=0$}\label{frogg}

There is a natural notion of generation
in the model.
We say that the first particle
in the origin
is of first generation
and that a particle
that turns red when it 
encounters a particle
of $k$th generation
is of $(k+1)$th generation.
(If a blue particle
moves on a site with
more than one red particles
then its generation
number is determined
by the lowest generation
number of the red particles.)
Now, to drive the red color
outside a ball an Euclidean ball
$B(0,r)$ in time $t$,
the first particle
initially in $z_1=0$
has to activate at some 
time $t_1$ and in some
site $z_2$ a second generation
particle,
this particle has to activate
at some time $t_1+t_2$
and in some site $z_3$
a third generation particle,\dots
and, for some $n$, an $n$th generation
will have to reach some site $z_{n+1}$
outside $B(0,r)$
at some time $t_1+\dots+t_n\leq t$.
Taking into account
the fact that more than one blue particle
can stand in a site reached
by a red particle
ad using Lemma \ref{dbb} we get,
for all $r$ and $t$,
\begin{equation}
P\left(
  \exists z \in {\cal R}(t),
  |z|>r
\right)
\leq Q(r,t)
\end{equation}
with 
\begin{eqnarray}
Q(r,t) &:=&
\sum_{n\geq 1}
\sum_{
  \stackrel{
    {\scriptstyle z_1,\dots, z_{n+1}}
   }{
    \stackrel{
      {\scriptstyle z_1=0}
    }{
      {\scriptstyle z_{n+1}\not\in B(0,r)}
    }
  }
} 
\int_{t_1+\dots + t_n\leq t}
\sum_{j_2,\dots,j_n\geq 0}
\prod_{k=2}^{n}
e^{-\rho}\frac{\rho^{j_k}}{j_k!}j_k
\nonumber\\
&&\quad
\prod_{k=1}^{n}
\left(
  \frac{\cst}{t_k^{d/2}}
  e^{-\frac{\cst|z_{k+1}-z_k|^2}{t_k}}
  \vee
  \cst e^{-\cst|z_{k+1}-z_k|}
\right)
dt_k
\label{dfnQ}
\end{eqnarray}
where here like in the sequel
we did not write,
to alleviate the notation,
that the integral is restricted
to positive variables only.

Permuting the last sum
with the product, making a spherical
change of variable
and using the triangular inequality
we get
\begin{equation}
Q(r,t) \leq
\sum_{n\geq 1}
\int_{
  \stackrel{
    {\scriptstyle r_1+\dots + r_n\geq r}
  }{
    {\scriptstyle t_1+\dots + t_n\leq t}
  }
}
\rho^{n-1}
\prod_{k=1}^{n}
q_1(r_k,t_k)\vee q_2(r_k) r_k^{d-1}
dr_kdt_k
\end{equation}
with
\begin{equation}
q_1(r_k,t_k) := 
  \frac{\cst}{t_k^{d/2}}
  e^{-\frac{\cst r_k^2}{t_k}}
\end{equation}
\begin{equation}
q_2(r_k,t_k) = q_2(r_k) := 
  \cst e^{-\cst r_k}
\end{equation}
Grouping together the different
terms according to the respective
values of $q_1$ and $q_2$
and using
\begin{equation}
\left(
  \begin{array}{c}
    \!\!\! n \!\!\!\\
    \!\!\! j \!\!\!\\
  \end{array}
\right)
\leq 2^n
\end{equation}
we get
\begin{eqnarray}
&& Q(r,t)
\nonumber\\
&&\quad\leq
\frac{1}{\rho}
\int_{
  \stackrel{
    {\scriptstyle R_1+R_2\geq r}
  }{
    {\scriptstyle T_1+T_2\leq t}
  }
}
\sum_{n\geq 1}
\sum_{j=0}^{n}
\left(
  \begin{array}{c}
    \!\!\! n \!\!\!\\
    \!\!\! j \!\!\!\\
  \end{array}
\right)
\left(
  \int_{
    \stackrel{
      {\scriptstyle r_1+\dots + r_j\geq R_1}
    }{
      {\scriptstyle t_1+\dots + t_j\leq T_1}
    }
  }
  \rho^{j}
  \prod_{k=1}^{j}
  q_1(r_k,t_k) r_k^{d-1}
  dr_kdt_k
\right)
\nonumber\\
&&\quad\quad\quad
\left(
  \int_{
    \stackrel{
      {\scriptstyle r_1+\dots + r_{n-j}\geq R_2}
    }{
      {\scriptstyle t_1+\dots + t_{n-j}\leq T_2}
    }
  }
  \rho^{n-j}
  \prod_{k=1}^{n-j}
  q_2(r_k) r_k^{d-1}
  dr_kdt_k
\right)
dR_1dR_2dT_1dT_2\\
&&\quad\leq
\frac{1}{\rho}
\int_{
  \stackrel{
    {\scriptstyle R_1+R_2\geq r}
  }{
    {\scriptstyle T_1+T_2\leq t}
  }
}
Q_1(R_1,T_1)Q_2(R_2,T_2)
dR_1dR_2dT_1dT_2
\end{eqnarray}
with for $j=1,2$
\begin{eqnarray}
Q_j(R_j,T_j) &:=&
\sum_{n\geq 1}
\int_{
  \stackrel{
    {\scriptstyle r_1+\dots + r_n\geq R_j}
  }{
    {\scriptstyle t_1+\dots + t_n\leq T_j}
  }
}
(2\rho)^{n}
\prod_{k=1}^{n}
q_j(r_k,t_k) r_k^{d-1}
dr_kdt_k
\end{eqnarray}
For any $R,T\geq 0$
we will estimate separately
$Q_1(R_1,T_1)$
and $Q_2(R_2,T_2)$.

We have
\begin{equation}
Q_1(R,T)\leq
\sum_{n\geq 1}
(\cst\rho)^n
\int_{
  \stackrel{
    {\scriptstyle r_1+\dots + r_n\geq R}
  }{
    {\scriptstyle t_1+\dots + t_n\leq T}
  }
}
\prod_{k=1}^{n}
e^{-\cst\frac{r_k^2}{t_k}}
\left(\frac{r_k}{\sqrt t_k}\right)^{d-1}
\frac{dr_k dt_k}{\sqrt t_k}
\end{equation}

Making a change of variable $x_k=\cst r_k^2/t_k$
and observing that, by the Cauchy-Schwartz inequality,
\begin{equation}
\left\{
  \begin{array}{l}
    \sum_k \sqrt t_k \sqrt x_k \geq \cst R\\
    \sum_k t_k \leq T
  \end{array}
\right.
\Rightarrow
\left\{
  \begin{array}{l}
    \sum_k x_k \geq \cst R^2/T\\
    \sum_k t_k \leq T
  \end{array}
\right.
\end{equation}
we get,
with $\Gamma$ the Euler function,
\begin{eqnarray}
Q_1(R,T)&\leq&
\sum_{n\geq 1}
(\cst\rho)^n
\int_{
  \stackrel{
    {\scriptstyle x_1+\dots + x_n\geq \cst R^2/T}
  }{
    {\scriptstyle t_1+\dots + t_n\leq T}
  }
}
\prod_{k=1}^{n}
e^{-x_k}
x_k^{\frac{d-1}{2}}
\frac{dx_k dt_k}{x_k^{1/2}}\\
&\leq&
\sum_{n\geq 1}
(\cst\rho)^n
\int_{
  \stackrel{
    {\scriptstyle x_1+\dots + x_n\geq \cst R^2/T}
  }{
    {\scriptstyle t_1+\dots + t_n\leq T}
  }
}
\prod_{k=1}^{n}
e^{-x_k}
x_k^{\frac{d}{2}-1}
\frac{dx_k dt_k}{\Gamma(d/2)}
\end{eqnarray}
Since the volume of the $n$-dimensional
simplex of side-length $T$ is $T^n/n!$
and the sum of independent variables
with a $\Gamma$ distribution
follows a $\Gamma$ law,
\begin{eqnarray}
Q_1(R,T)
&\leq&
\sum_{n\geq 1}
\frac{(\cst\rho T)^n}{n!}
\int_{x\geq \cst R^2/T}
e^{-x}
x^{n\frac{d}{2}-1}
\frac{dx}{\Gamma\left(n\frac{d}{2}\right)}\\
&\leq&
\sum_{n\geq 1}
\frac{(\cst\rho T)^n}{n!}
P\left(N'\leq\left\lceil\frac{nd}{2}\right\rceil\right)\\
&\leq&
e^{\cst\rho T}
P\left(N'\leq\cst N\right)
\end{eqnarray}
where $N$ and $N'$
are independent Poissonian variables
of mean $\cst R^2/T$ respectively.
Now, for any large enough $\gamma$,
if $R\geq\gamma\sqrt\rho T$,
then by (\ref{ldiii})
\begin{eqnarray}
Q_1(R,T)
&\leq&
e^{\cst\rho T}
P\left(
  \frac{N}{E[N]}
  \geq
  \cst\frac{R^2/T}{\rho T}\frac{N'}{E[N']}
\right)\\
&\leq&
e^{\cst\rho T}
e^{-\cst\frac{R^2}{T}}\\
&\leq&
e^{\cst\rho T}
e^{-\cst\sqrt\rho R}
\end{eqnarray}
so that, for any large enough $\gamma$,
\begin{eqnarray}
Q_1(R,T)
&\leq&
e^{\cst\rho T}
\exp\left\{
  -\cst\sqrt\rho R
  \one_{[\gamma\sqrt\rho T,+\infty[}(R)
\right\}
\end{eqnarray}

Turning to $Q_2(R,T)$
we have
\begin{eqnarray}
Q_2(R,T)
&\leq &
\sum_{n\geq 1}
(\cst\rho)^n
\int_{
  \stackrel{
    {\scriptstyle r_1+\dots + r_n\geq R}
  }{
    {\scriptstyle t_1+\dots + t_n\leq T}
  }
}
\prod_{k=1}^{n}
e^{-\cst r_k}
r_k^{d-1}
dr_k dt_k\\
&\leq&
\sum_{n\geq 1}
(\cst\rho)^n
\int_{
  \stackrel{
    {\scriptstyle x_1+\dots + x_n\geq \cst R}
  }{
    {\scriptstyle t_1+\dots + t_n\leq T}
  }
}
\prod_{k=1}^{n}
e^{-x_k}
x_k^{{d-1}}
dx_k dt_k\\
&\leq&
\sum_{n\geq 1}
\frac{(\cst\rho T)^n}{n!}
\int_{x\geq \cst R}
e^{-x}
x^{nd-1}
\frac{dx}{\Gamma(nd)}\\
&\leq&
e^{\cst\rho T}
P\left(N'\leq\cst N\right)
\end{eqnarray}
where $N$ and $N'$
are independent Poissonian variables
of mean $\cst R$ respectively.
Then, for any large enough $\gamma$,
if $R\geq\gamma\rho T$,
we get by (\ref{ldiii})
\begin{eqnarray}
Q_2(R,T)
&\leq&
e^{\cst\rho T}
P\left(
  \frac{N}{E[N]}
  \geq
  \frac{\cst R}{\rho T}\frac{N'}{E[N']}
\right)\\
&\leq&
e^{\cst\rho T}
e^{-\cst R}
\end{eqnarray}
so that, for any large enough $\gamma$,
\begin{eqnarray}
Q_2(R,T)
&\leq&
e^{\cst\rho T}
\exp\left\{
  -\cst R
  \one_{[\gamma\rho T,+\infty[}(R)
\right\}
\end{eqnarray}

Turning back to Q(r,t),
we get, for any large enough $\gamma$,
\begin{eqnarray}
&& Q(r,t)
\nonumber\\
&&\quad\leq\quad
\frac{1}{\rho}
\int_{
  \stackrel{
    {\scriptstyle R_1+R_2\geq r}
  }{
    {\scriptstyle T_1+T_2\leq t}
  }
}
e^{\cst\rho (T_1+T_2)}
\nonumber\\
&&\quad\quad\quad\quad
\exp\left\{
  -\cst\left(
    \sqrt\rho R_1
    \one_{[\gamma\sqrt\rho T_1,+\infty[}(R_1)
    +
    R_2
    \one_{[\gamma\rho T_2,+\infty[}(R_2)
  \right)
\right\}
\nonumber\\
&&\quad\quad\quad\quad\quad
dR_1dR_2dT_1dT_2\\
&&\quad\leq\quad
\frac{1}{\rho}
\int_{
  \stackrel{
    {\scriptstyle R_1+R_2\geq r}
  }{
    {\scriptstyle T_1+T_2\leq t}
  }
}
e^{\cst\rho t}
\nonumber\\
&&\quad\quad\quad\quad
\exp\left\{
  -\cst\left(
    \sqrt\rho R_1
    \one_{[\gamma\rho T_1,+\infty[}(\sqrt\rho R_1)
    +
    R_2
    \one_{[\gamma\rho T_2,+\infty[}(R_2)
  \right)
\right\}
\nonumber\\
&&\quad\quad\quad\quad\quad
dR_1dR_2dT_1dT_2
\end{eqnarray}

Now if $\rho\leq 1$, then
\begin{equation}
R_2
\one_{[\gamma\rho T_2,+\infty[}(R_2)
\geq
\sqrt\rho R_2
\one_{[\gamma\rho T_2,+\infty[}(\sqrt\rho R_2)
\end{equation}
and if $\rho\geq 1$, then
\begin{equation}
\sqrt\rho R_1
\one_{[\gamma\rho T_1,+\infty[}(\sqrt\rho R_1)
\geq
R_1
\one_{[\gamma\rho T_1,+\infty[}(R_1)
\end{equation}

As a consequence,
with
\begin{equation}
\bar\rho:= \max(\rho,\sqrt\rho)
\mbox{ and }
X_i=\frac{\rho}{\bar\rho}R_i,
\quad
i=1,2
\end{equation}
we have
\begin{eqnarray}
&&Q(r,t)
\nonumber\\
&&\quad\leq\quad
\frac{\bar\rho^2}{\rho.\rho^2}
\int_{
  \stackrel{
    {\scriptstyle X_1+X_2\geq \rho r/\bar\rho}
  }{
    {\scriptstyle T_1+T_2\leq t}
  }
}
e^{\cst\rho t}
\nonumber\\
&&\quad\quad\quad\quad
\exp\left\{
  -\cst\left(
    X_1
    \one_{[\gamma\rho T_1,+\infty[}(X_1)
    +
    X_2
    \one_{[\gamma\rho T_2,+\infty[}(X_2)
  \right)
\right\}
\nonumber\\
&&\quad\quad\quad\quad\quad
dX_1dX_2dT_1dT_2\\
&&\quad\leq\quad
\frac{e^{\cst\rho t}}{\rho^2}
\int_{
  \stackrel{
    {\scriptstyle X_1+X_2\geq \rho r/\bar\rho}
  }{
    {\scriptstyle T_1+T_2\leq t}
  }
}
e^{-\cst\left(
  X_1+X_2
  -\gamma\rho(T_1+T_2)
\right)}
\nonumber\\
&&\quad\quad\quad\quad
dX_1dX_2dT_1dT_2
\end{eqnarray}

If $r\geq 2\gamma\bar\rho t$, i.e.,
\begin{equation}
\frac{\rho r}{2\bar\rho}
\geq 
\gamma\rho t 
\end{equation}
then
\begin{eqnarray}
Q(r,t)
&\leq&
\frac{e^{\cst\rho t}}{\rho^2}
\int_{
  \stackrel{
    {\scriptstyle X_1+X_2\geq \rho r/\bar\rho}
  }{
    {\scriptstyle T_1+T_2\leq t}
  }
}
e^{-\cst\frac{X_1+X_2}{2}} 
dX_1dX_2dT_1dT_2\\
&\leq&
\cst
\frac{e^{\cst\rho t}}{\rho^2}
t^2
\frac{\rho r}{\bar\rho}
e^{-\cst\frac{\rho r}{2\bar\rho}}\\
&\leq&
\cst
\frac{e^{\cst\rho t}}{\rho^4}
e^{-\cst\frac{\rho r}{2\bar\rho}}\\
&\leq&
\frac{\cst}{\rho^4}
e^{\cst\rho t}
e^{-\cst\gamma\rho t} 
\end{eqnarray}
and, with a large enough $\gamma$, 
we get
\begin{eqnarray}
Q(r,t)
&\leq&
\frac{\cst}{\rho^4}
e^{-\cst\rho t}
\end{eqnarray}
\qed

\subsection{Proof of Proposition \ref{pf1}}
\label{frog1}

In the previous proof
we could have use,
instead of the estimates from
Lemma \ref{dbb} on 
$P_0(\zeta(t)=z)dt$,
an estimate on
\begin{equation}
dP_0\left(
  \tau_z(\zeta)\leq t
\right)
=
P_0\left(
  \tau_z(\zeta)\in [t,t+dt]
\right)
\end{equation}
with
\begin{equation}
\tau_z(\zeta):=
\inf\left\{
  t\geq 0 :\:
\zeta(t)=z
\right\}
\end{equation}
While in dimension
$d\geq 2$
the two quantities are quite close,
in dimension $d=1$
they are substantially different.
In addition,
using $\tau_z(\zeta)$
in dimension 1
allows for a simpler proof
of a stronger result when $\rho$
is small enough.
Indeed, for all $r$ and $t$,
\begin{eqnarray}
&&P\left(
  \exists z\in {\cal R}(t),
  |z|>r
\right)
\nonumber\\
&&\quad\leq\quad
\sum_{n\geq 1}
\rho^{n-1}
\sum_{r_1+\dots r_n \geq r}
\int_{t_1+\dots + t_n\leq t}
\prod_{k=1}^{n}
dP_0\left(
  \tau_{r_k}(\zeta)\leq t_k
\right)\\
&&\quad\leq\quad
\sum_{n\geq 1}
\rho^{n-1}
\sum_{R\geq r}
\sum_{r_1+\dots r_n = R}
P_0\left(
  \tau_{R}(\zeta)\leq t
\right)
\end{eqnarray}
Then, by the reflexion principle
and Lemma \ref{dbb}
\begin{eqnarray}
P\left(
  \exists z\in {\cal R}(t),
  |z|>r
\right)
&\leq&
\frac{\cst}{\rho}
\sum_{R\geq r}
\sum_{n\geq 1}
\frac{\rho^{n}R^n}{n!}
\left(
  e^{-\cst R^2/T}\vee e^{-\cst R}
\right)\\
&\leq&
\frac{\cst}{\rho}
\sum_{R\geq r}
e^{\rho R}\left(
  e^{-\cst R^2/T}\vee e^{-\cst R}
\right)
\end{eqnarray}
Now if $r\geq\gamma\rho t$
for some large enough $\gamma$
we get,
for $\rho$ small enough,
\begin{eqnarray}
P\left(
  \exists z\in {\cal R}(t),
  |z|>r
\right)
&\leq&
\frac{\cst}{\rho}
\sum_{R\geq r}
e^{-\cst\rho R}\\
&\leq&
\frac{\cst}{\rho^2}
e^{-\cst\rho r}\\
&\leq&
\frac{\cst}{\rho^2}
e^{-\cst\rho^2 t}
\end{eqnarray}
This proves Proposition
\ref{pf1}
for small $\rho$'s.
When $\rho$
is bounded away
from~0,
Proposition \ref{pf1}
is just a consequence of Theorem \ref{thm1}
for the frog model.

\section{RB and RBK processes}
\label{gc}

\subsection{Proof of theorem \ref{thm1}}

We can proceed like in the case of the frog model
except for the fact that a particle
does not anymore turns red
at the same point where it started.
We have then to sum on the possible starting points.
With the notation
\begin{equation}
s_k=t_1+\dots +t_{k-1}, \quad k\geq 2
\end{equation}
and for any $i\geq 1$
we have
\begin{eqnarray}
&& P\left(
  \exists z \in {\cal R}(t),
  |z|>r
\right)\\
&& \quad\leq
\sum_{n\geq 1}
\sum_{
  \stackrel{
    {\scriptstyle z_1,\dots, z_{n+1}}
   }{
    \stackrel{
      {\scriptstyle z_1=0}
    }{
      {\scriptstyle z_{n+1}\not\in B(0,r)}
    }
  }
}
\int_{t_1+\dots + t_n\leq t}
\sum_{
  \stackrel{
    {\scriptstyle z'_2,\dots, z'_{n}}
  }{ 
    {\scriptstyle j_2,\dots,j_n\geq 0}
  }
}
\prod_{k=2}^{n}
e^{-\rho}\frac{\rho^{j_k}}{j_k!}j_k
P(z'_k+Z_i^B(s_k)=z_k)
\nonumber\\
&&\quad\quad\quad
\prod_{k=1}^{n}
\left(
  \frac{\cst}{t_k^{d/2}}
  e^{-\frac{\cst|z_{k+1}-z_k|^2}{t_k}}
  \vee
  \cst e^{-\cst|z_{k+1}-z_k|}
\right)
dt_k
\end{eqnarray}

Now permuting
the last sum
with the product
and using (\ref{bs}) we get
\begin{equation}
P\left(
  \exists z \in {\cal R}(t),
  |z|>r
\right)
\leq
Q(r,t)
\end{equation}
with $Q(r,t)$
defined in (\ref{dfnQ})
and estimated in the previous section.
\qed

\noindent
{\bf Remark:}
Unfortunately the proof of Proposition \ref{pf1}
cannot be extended so simply
to the general case,
even if
we restrict ourselves to KS processes.
To do so we would have to link the differential
\begin{equation}
dP_0\left(
  \tau_{z_R}(\zeta_R)\leq t
\right)
=
P_0\left(
  \tau_{z_R}(\zeta_R)\in [t,t+dt]
\right)
\end{equation}
with the sum
\begin{equation}
\sum_{z_B>0}
P_{(0,z_B)}\left(
  \tau_0(\zeta_B-\zeta_R)\in [t,t+dt],
  \zeta_R(t)=z_R
\right)
\end{equation}
with $\zeta_R$ and $\zeta_B$
independent continuous time random walks
with jump rates $D_R=1$
and $D_B>0$.
In the case $D_B=1$ this can be done
using the independence
between $\zeta_B-\zeta_R$ and
$\zeta_B+\zeta_R$.
In the case $D_B \neq 1$
we can only use an ``asymptotic independence''
between $\zeta_B-\zeta_R$ and
$\zeta_B+D_B\zeta_R$.
In both cases this is a quite technical task:
we will not go in this paper
beyond the result
for the frog model.

\subsection{Proof of Theorem \ref{thm2}}

We can adapt the proof for the frog model
using the QRW property and the last observations
of section \ref{qrw}:
\begin{eqnarray}
&&
P\left(
  \exists z \in {\cal R}(T),
  |z|>R,
  T>\Tau_{\alpha,\lambda}
\right)
\nonumber\\
&& \quad\leq\quad 
\sum_{n = 1}^{\lceil\lambda l T/T_\alpha\rceil}
\sum_{
  \stackrel{
    {\scriptstyle z_1,\dots, z_{n+1} \in \Lambda_\beta}
   }{
    \stackrel{
      {\scriptstyle z_1=0}
    }{
      {\scriptstyle z_{n+1}\not\in B(0,R)}
    }
  }
} 
\int_{t_1+\dots + t_n\leq T}
\prod_{k=1}^{n}
\cst\lambda^3 l^2
\nonumber\\
&&\quad\quad\quad
\left(
  \frac{\cst}{t_k^{d/2}}
  e^{-\frac{\cst|z_{k+1}-z_k|^2}{t_k}}
  \vee
  \cst e^{-\cst|z_{k+1}-z_k|}
\right)
dt_k
+SES
\end{eqnarray}
In this formula the first sum
is limited to $\lceil\lambda l T/T_\alpha\rceil$
since $T$ is at most exponential in $\beta$
and in each interval of length $T_\alpha$,
with probability $1-SES$,
interactions are limited to clouds
that contains $\lambda$ particles
at most and particles
are coupled with random walks
with $l$ pauses at most.
The factor  $l^2$ is due to the fact
that,
with probability $1-SES$,
in each pause interval the distance
between a particle and its associated
random walk with pauses increases
of $l$ at most,
one factor $\lambda$
is due to the fact that $\lambda$
red particles at most can leave
a given cluster before $\Tau_{\alpha,\lambda}$
and the last factor $\lambda^2$
is due to the fact that at each time
$t<\Tau_{\alpha,\lambda}$
a given particle
can turn red other particles
inside a radius $\lambda$
at most.

Then we can repeat the calculation of section
\ref{frogg}
with two main differences.
On the one hand we do not have
anymore the factor $\rho^{n-1}$ in our sum,
on the other hand this sum
is limited to 
$\lceil\lambda l T/T_\alpha\rceil$.
Instead of (\ref{ldiii})
we use then (\ref{ldii}) repeatedly.
For example
defining $Q_1$
and $Q_2$ in an analogous way
and observing that for any $\delta >0$,
$\lambda$ and $l$ are smaller
$e^{\delta\beta}$
for $\beta$ large enough, we have now
\begin{equation}
Q_1(R,T)\leq
\sum_{n=1}^{\lceil e^{\delta\beta}\rho T\rceil}
\frac{
  \left(
     e^{\delta\beta}T
  \right)
}{
  n!
}
P\left(
  N'\leq
  e^{\delta\beta}\rho T
\right)
+ SES
\end{equation}
with $N'$
a Poisson variable of mean
$\cst R^2/T$.
For any $\delta_1>\delta/2$,
if $R\geq e^{\delta_1\beta}\sqrt\rho T$
the last probability
can be estimated from above by
\begin{equation}
P\left(
  N'\leq
  e^{\delta\beta}\rho T
\right)
\leq
\exp\left\{
  -\cst\rho t
\right\}
+ SES
\end{equation}
while the last
sum can be estimated from above by
\begin{equation}
\sum_{n=1}^{\lceil e^{\delta\beta}\rho T\rceil}
\frac{
  \left(
     e^{\delta\beta}T
  \right)
}{
  n!
}
\leq
\exp\{e^{\delta\beta}T\}
P\left(
  N\leq
  \cst e^{\delta\beta}\rho T
\right)
+ SES
\end{equation}
with $N$
a Poisson variable
of mean $e^{\delta\beta}\rho T$,
so that,
by (\ref{ldii}),
\begin{equation}
\sum_{n=1}^{\lceil e^{\delta\beta}\rho T\rceil}
\frac{
  \left(
     e^{\delta\beta}T
  \right)
}{
  n!
}
\leq
\exp\{e^{2\delta\beta}\rho T\}
+ SES
\end{equation}
Putting everything together
we get,
for any $R$, $T$,
\begin{equation}
Q_1(R,T)\leq
\exp\{e^{2\delta\beta}\rho T\}
\exp\left\{
  -\cst\sqrt\rho R
  \one_{[e^{\delta_1\beta}\sqrt\rho T, +\infty[}
  (R)
\right\} + SES
\end{equation}
We can estimate
$Q_2$ in the same way
and the rest of the calculation
goes like in section \ref{frogg}.
\qed

\subsection{Proof of Theorem \ref{thm3}}

Given $\Lambda^{(1)}$ and 
$\Lambda^{(2)}$
with
\begin{equation}
d(\Lambda^{(1)},\Lambda^{(2)})
> e^{\delta\beta}\max(\sqrt T, \sqrt\rho T)
\end{equation}
we define a new coloring process.
With
\begin{equation}
B:=\Lambda^{(1)}\cup\Lambda^{(2)}
\end{equation}
and 
\begin{equation}
W:=\left\{
  z\in\Lambda_\beta :\:
  \inf_{b\in B}
  |z-b|
  > e^{-\delta\beta/2}d(\Lambda^{(1)},\Lambda^{(2)})
\right\}
\end{equation}
we say that all the particles
that start from $B$ are black,
all the particles that start from $W$
are white and all the particles that start
from $(B\cup W)^c$
do not have any color
at time $t=0$
Then, for $t>0$,
black particles keep their black color,
white particles keep their white color,
non-colored particles that enter $B$ turn black,
non-colored particles that enter $W$ turn white,
and 
non-colored particles that share some cluster 
with a colored particle
turn black or white
choosing randomly a colored particle
inside the cluster and taking
the same color.
We can define a black zone
and a white zone
like we defined 
the red zone.
As a consequence of corollary
\ref{kcor}, with probability
$1-SES$,
the black and white zones
will not intersect
up to time
$T\wedge\Tau_{\alpha,\lambda}$
and
we will never see
black and white particles
in a same cluster
up to time
$T\wedge\Tau_{\alpha,\lambda}$.

Now we couple in the more natural way
the previous process,
with a process that starts
from the same initial configuration,
uses the same marks and clocks
for the particles
and evolves in the same way
except for the fact
that 
each particle in $W$
or that enters in $W$
disappears.
For this process the restrictions of the dynamics
to $\Lambda^{(1)}$ and 
$\Lambda^{(2)}$
are clearly independent
and the previous observation
shows that,
with probability
$1-SES$,
these restrictions 
for the two processes
coincide
up to time $T\wedge\Tau_{\alpha,\lambda}$.
This proves the theorem.
\qed

\section*{Acknowledgments}

We thank Eurandom
for its hospitality
and the Grefi-Mefi
for partial support
as well as for the organization
of its 2008 Workshop
that stimulated lot of this work.

We thank Francesco Manzo for his idea
to estimate the propagation velocity
looking at the number of particles.
This is important
because it founds
the heuristic that allows
us to argue we proved ``good bounds''.
We thank Amine Asselah
for his help in correcting
a previous and wrong version
of Lemma \ref{dbb}. 
We thank Beatrice Nardi for her hospitality.
We thank Pietro Glasmacher for his enthusiasm
during our work sessions
and Patrick Glasmacher for his support
during the same work sessions.

\end{document}